\documentclass[draft,12pt]{article}

\usepackage[latin1]{inputenc}
\usepackage{amsmath,amssymb}
\usepackage{latexsym}

\newtheorem{theorem}{Theorem}[section]
\newtheorem{corollary}[theorem]{Corollary}
\newtheorem{lemma}[theorem]{Lemma}
\newtheorem{proposition}[theorem]{Proposition}
\newtheorem{definition}[theorem]{Definition}
\newtheorem{assumption}[theorem]{Assumption}
\newtheorem{remark}[theorem]{Remark}

\numberwithin{equation}{section}

\def\proof{{\medskip\noindent {\bf Proof. }}}
\def\longproof#1{{\medskip\noindent {\bf Proof #1.}}}
\def\qed{{\hfill $\square$ \bigskip}}

\def\square{{\vcenter{\vbox{\hrule height.3pt
        \hbox{\vrule width.3pt height5pt \kern5pt
           \vrule width.3pt}
        \hrule height.3pt}}}}

  \def\sC {{\cal C}}
 \def\sE {{\cal E}} \def\sF {{\cal F}}
\def\sG {{\cal G}}  
  \def\sL {{\cal L}}

\def\sV {{\cal V}}

\def\wt{\widetilde}

\def\E{{\mathbb E}}
\def\P{{\mathbb P}}
\def\norm#1{{\Vert #1 \Vert}}
\def\del{{\partial}}
\def\lam{{\lambda}}
\def\angel#1{{\langle #1 \rangle}}
\def\bee{\begin{equation}}
\def\bet{\begin{theorem}}
\def\bep{\begin{proposition}}
\def\bel{\begin{lemma}}
\def\bec{\begin{corollary}}
\def\bed{\begin{definition}}
\def\ber{\begin{remark}}
\def\eee{\end{equation}}
\def\eet{\end{theorem}}
\def\eep{\end{proposition}}
\def\eel{\end{lemma}}
\def\eec{\end{corollary}}
\def\eed{\end{definition}}
\def\eer{\end{remark}}

\def\R{{\mathbb R}}
\def\E{{{\mathbb E}\,}}
\def\P{{\mathbb P}}

\def\Z{{\mathbb Z}}

\def\lam{{\lambda}}
\def\th{{\theta}}
\def\al{{\alpha}}
\def\grad{{\nabla}}
\def\proof{{\medskip\noindent {\bf Proof. }}}
\def\longproof#1{{\medskip\noindent {\bf Proof #1.}}}
\def\qed{{\hfill $\square$ \bigskip}}

\def\vp{\varphi}

\def\angel#1{{\langle#1\rangle}}
\def\norm#1{\Vert #1 \Vert}

 \def\qq {\qquad}
\def\del{{\partial}}
\def\wt{\widetilde}

\def\ni{\noindent }
\def\ms{\medskip}

\def\square{{\vcenter{\vbox{\hrule height.3pt
        \hbox{\vrule width.3pt height5pt \kern5pt
           \vrule width.3pt}
        \hrule height.3pt}}}}

\def\tfrac#1#2{{\textstyle {\frac{#1}{#2}}}}

\def\tlint{{- \kern-0.85em \int \kern-0.2em}}  
\def\dlint{{- \kern-1.05em \int \kern-0.4em}}  

  \def\sC {{\cal C}}
 \def\sE {{\cal E}} \def\sF {{\cal F}}
\def\sG {{\cal G}}  
  \def\sL {{\cal L}}

\def\sV {{\cal V}}

\def\nn{{\nonumber}}

\begin{document}

\title{A stability theorem for  elliptic Harnack inequalities}
\author{Richard F. Bass\footnote{Research partially supported by NSF grant
DMS-0901505.}}

\date{\today}

\maketitle

\begin{abstract}  
\noindent {\it Abstract:} 
We prove a stability theorem for the elliptic Harnack inequality:
if two weighted graphs are equivalent,
then the elliptic Harnack inequality holds for harmonic functions with
respect to one of the graphs if and only if it holds for harmonic functions
with respect to the other graph.
As part of the proof, we give a characterization of the elliptic Harnack inequality.

\vskip.2cm
\noindent \emph{Subject Classification: Primary 31B05; Secondary 31E05, 60J27}   
\end{abstract}

\section{Introduction}\label{S:intro}

A justly famous theorem of Moser \cite{Moser1} says that if
 $\sL$ is the  uniformly elliptic operator in divergence form
given by
$$\sL f(x)=\sum_{i,j=1}^d \frac{\del }{\del x_i}\Big(a_{ij}(\cdot)
\frac{\del f}{\del x_j}(\cdot)\Big)(x)$$
acting on functions on  $\R^d$, 
where the $a_{ij}$ are also bounded and measurable, then an
elliptic Harnack inequality (EHI) holds for functions that are non-negative
and harmonic with respect to $\sL$ in a domain. This is one of the more
important theorems in the study of elliptic and parabolic partial differential
equations,
and is used, for example, in deriving \emph{a priori} regularity results
for harmonic functions and for heat kernels.

The operator $\sL$ is associated
with the Dirichlet form
$$\sE_\sL(f,f)=\int_{\R^d} \sum_{i,j=1}^d a_{ij}(x) \frac{\del f}{\del x_i}(x)
\frac{\del f}{\del x_j}(x)\, dx.$$
If the $a_{ij}$'s are bounded and the matrices $a(x)=(a_{ij}(x))$
are uniformly positive definite, then $\sE_\sL$ is comparable to 
$\sE_\Delta$, where
$$\sE_\Delta(f,f)=\int_{\R^d} |\grad f(x)|^2\, dx,$$
which is the Dirichlet form corresponding to the Laplacian. Thus 
one could rephrase Moser's theorem as saying that whenever the 
Dirichlet form corresponding to an operator $\sL$ is comparable to
the Dirichlet form corresponding to the Laplacian, then the EHI
 holds for non-negative functions that are harmonic
with respect to $\sL$ in a domain. 

We can view Moser's theorem  as a stability theorem for the EHI. The purpose of
this paper is to generalize this stability property to very general 
state spaces. We show that provided 
some mild regularity holds, then whenever
two Dirichlet forms $\sE_1$ and $\sE_2$ are comparable with corresponding 
operators $\sL_1$ and $\sL_2$,  the EHI
holds for non-negative harmonic functions with respect
to $\sL_1$  if and only if the EHI  holds
for non-negative harmonic functions with respect
to $\sL_2$.

We also provide a characterization of the EHI.
Provided the regularity  holds, this characterization can be
considered as necessary and sufficient conditions for the EHI.

It is interesting to compare the EHI with the parabolic Harnack
inequality (PHI). The PHI, first proved by Moser in \cite{Moser2}
(see \cite{Fabes-Stroock} for a very different proof), is a Harnack
inequality for non-negative solutions to 
$$\frac{\del u}{\del t}(x,t)=\sL u(x,t)$$
in a domain. Necessary and sufficient conditions are known for the PHI
in quite general state spaces. If the state space is regular enough to
have a large class of nice cut-off functions, then 
Grigor'yan \cite{Grigoryan}  and
Saloff-Coste \cite{Saloff-Coste}
independently proved that the PHI
holds if and only if both volume doubling  and a Poincar\'e inequality 
hold. This was extended to the case where such nice cut-off functions
need not exist in \cite{BB} and \cite{BBK}. The latter papers allow
state spaces that have fractal structure or that have large numbers
of obstructions.   

If the PHI holds, then the EHI holds; this is quite easy to see.
The converse 
is false. In \cite{BB1} an example was given where EHI holds, but
the PHI in the usual form does not (that is, with scaling factor $r^2$). Delmotte \cite{Delmotte} constructed
an example where the EHI holds, but volume doubling does not, and 
consequently the PHI cannot hold in any form. See \cite{Hebisch-Saloff-Coste}
for more on the relationship between the EHI and PHI.
It has been an open problem for quite some time to find a characterization
of the EHI comparable to the one for the PHI.

In this paper we primarily 
look at infinite graphs rather than continuous state
spaces. All the key ideas are present in the infinite graph case and
we avoid some unpleasant technicalities. It is straightforward to
extend our results to metric measure Dirichlet spaces in a manner
very similar to how \cite{BBK} extended \cite{BB}; see Section \ref{S:further}.  

We consider infinite graphs where between any two adjacent vertices
$x$ and $y$ there is given a conductance $C_{xy}$. If $x$ and $y$ are adjacent,
we write $x\sim y$. Setting
$$\mu_x=\sum_{z\sim x} C_{xz},$$
we can construct a continuous time Markov chain $X$ with the graph
as the state space. When $X$ is at $x$, it waits an independent
exponential length of time with parameter $\mu_x$  and 
then jumps to an adjacent vertex.
It chooses a neighboring vertex $y$ with probability $C_{xy}/\mu_x$.   
We write $\sL$ for the infinitesimal generator of $X$. A function
$h$ is harmonic with respect to $\sL$ in a domain $D$ if
$$h(x)=\sum_{y\sim x} h(y) C_{xy}, \qq x\in D.$$
Let $B(x,r)$ denote the ball of radius $r$ about $x$.
The elliptic Harnack inequality states that there exists a constant
$c$ not depending on $x_0$ or $r$ such that 
if $h$ is non-negative and harmonic in $B(x_0,2r)$, then 
$$h(x)\le ch(y), \qq x,y\in B(x_0,r).$$

We do require some mild regularity. 
 For example, one of our assumptions is that volume doubling holds.
Whereas the PHI implies volume doubling, the example of Delmotte
\cite{Delmotte} shows that the EHI can hold even though volume doubling
does not. Since every known approach to proving an EHI uses volume
doubling in an essential way, the problem of finding  necessary
and sufficient conditions for the EHI to hold without assuming any
regularity looks very hard.

For most of this paper we consider the case where the process $X$ is  transient. That is, $d(X_t,x)\to \infty$
almost surely as $t\to \infty$ for every point $x$, where $d(\cdot, \cdot)$
is the graph distance. This, for example, allows us to define capacities.
The general case, which is slightly more complicated to state, is
given in Section \ref{S:further}.

Let 
$V(x,r)$ be the volume of $B(x,r)$ with respect to the measure 
$\mu(A)=\sum_{x\in A} \mu_x$. Let $C(x,r)$ be the capacity of 
$B(x,r)$ (a definition is given  in the next section).
Finally define $E(x,r)=V(x,r)/C(x,r)$. It will turn out that
$E(x,r)$ is comparable to the expected time that the process spends 
in $B(x,r)$ when started at $x$.

The novel feature of this paper is to introduce the
adjusted Poincar\'e inequality (API):
$$\sum_{y\in B(x,r)} |f(y)-f_{B(x,r)}|^2\, \mu_y\leq cE(x,r)\sE_{B(x,c'r)}(f,f).$$
Here $c'>1$, $f_A$ is the average value of $f$ on the set $A$ with respect to
the measure $\mu$, and $\sE_A$ is the Dirichlet form restricted
to the set $A$. Note that in the usual Poincar\'e inequality, $E(x,r)$
is replaced by $r^\beta$ for $\beta$ equal to some constant, most often,
$\beta=2$.

We will also use another inequality, which we call the  cut-off inequality (COI). This is closely related
to the cut-off Sobolev inequality of \cite{BB}. 

Our first main theorem is that if transience and regularity hold, then
the EHI holds if and only both the COI and API hold. This immediately implies
our second theorem, the stability result, which says that if  
transience and the regularity hold and the EHI holds for a weighted graph, then the EHI holds
for every equivalent weighted graph.
These results are new even when sufficiently many nice cut-off functions
exist.

In the next section we give a precise statement of our results. In 
Section \ref{S:prelim} we introduce the cable process and also prepare
some preliminary results. Section \ref{S:conseq} proves some estimates
that can be obtained from the EHI. We prove that the EHI implies the API
in Section \ref{S:api}, and prove our main theorems in
Section \ref{S:proofs}.
In Section \ref{S:further} we consider the general case (where $X$ is not
necessarily transient). In that section we also consider there extensions to the situation
where the state space is a metric measure space rather than a graph.

\section{Statement of results}\label{S:sor}

We use the letter $c$ with subscripts to denote finite positive
constants whose exact values are unimportant and may change from
place to place.

Let $\sG$ be an infinite connected graph consisting of  vertices $\sV$ together
with a collection of edges. We write
$x\sim y$ if $x$ and $y$ are vertices connected by an edge. We suppose each
vertex belongs to at most finitely many edges. 
For each pair $x,y\in \sV$  we define a conductance
$C_{xy}\ge 0$ 
such that $C_{xy}=C_{yx}$ and also $C_{xy}=0$ unless
$x\sim y$. The graph $\sG$ together with the conductances $\{C_{xy}\}$
is called a weighted graph.

Let $\mu_x=\sum_{y} C_{xy}$,
and define a measure $\mu$ on $\sV$ by $\mu(A)=\sum_{x\in A} \mu_x$. 
We let $d(x,y)$ be the usual graph distance on $\sG$ and set
$$B(x,r)=\{y: d(x,y)<r\}, \qq V(x,y)=\mu(B(x,r)).$$
We assume throughout this paper that there exists a constant $c_1$
such that 
\bee\label{p0}
0<\mu_x\le c_1, \qq x\in \sV.
\eee

For $f\in L^2(\sV, \mu)$, define
$$\sE_\sG(f,f)=\tfrac12 \sum_{x\sim y} [f(y)-f(x)]^2 C_{xy}$$
and 
$$\sF_\sG=\{f\in L^2(\sV,\mu): \sE_\sG(f,f)<\infty\}.$$
It is well known (see \cite{FOT}) that $(\sE_\sG, \sF_\sG)$ is
a regular Dirichlet form associated with a strong Markov process $(X_t, \P^x)$.
The process $X$ is a continuous time Markov chain on $\sV$ which can be described
as follows. When $X$ is at a vertex $x$, it waits there an independent
exponential length of time with parameter $\mu_x$ and then jumps to a
neighboring vertex. It chooses the neighboring vertex $y$ to jump to with 
probability $C_{xy}/\mu_x$. The infinitesimal generator of $X$ is given by
$$\sL_\sG f(x)=\sum_{x\sim y} [f(y)-f(x)]C_{xy}.$$

Except for Section \ref{S:further}  we make  a transience assumption.

\begin{assumption}\label{A1}{\rm
$(\sE_\sG, \sF_\sG)$ is transient in the sense of \cite[Sect.~1.5]{FOT}.
}
\end{assumption}

An equivalent formulation in our context is that 
$$\lim_{t\to \infty} d(X_t,x)\to \infty$$
with probability one for each starting point and each $x\in \sV$.

Let $$C(x,r)= \inf\{\sE_\sG(f,f): f\in \sF_\sG, f|_{B(x,r)}=1\}$$
be the capacity of $B(x,r)$. This exists and is finite because
$(\sE_\sG, \sF_\sG)$ is transient; see \cite[Sect.~2.1]{FOT}. Define
\bee\label{sor-E1}
E(x,r)=\frac{V(x,r)}{C(x,r)}.
\eee
We will see later that $E(x,r)$ is comparable to the expected occupation
time of $B(x,r)$ by $X_t$ when started at $x$.

Our second main assumption concerns regularity.

\begin{assumption}\label{A2}{\rm
There exist $c_1>0$ and $\rho\in (0,1)$ such that
the following three inequalities hold.

\ni Volume doubling holds:
\bee\label{A2-VD}
V(x,2r)\le c_1 V(x,r), \qq x\in \sV, r\ge 1.
\eee
Capacity growth holds:
\bee\label{A2-CG}
C(x,r)\le \rho C(x,2r), \qq x\in \sV, r\ge 1.
\eee
Expected occupation time growth holds:
\bee\label{A2-EG}
E(x,r)\le \rho E(x,2r), \qq x\in \sV, r\ge 1.
\eee
}
\end{assumption}

Finally we need a geometric condition.

\begin{assumption}\label{A3}{\rm
There exists $M$ not depending on $x$ or $r$ such that
the boundary of $B(x,r)$ can be covered by at most $M$ balls of radius
$r/8$ provided $r\ge 1$.
}
\end{assumption}

Regarding our assumptions, we make these remarks.

\begin{remark}\label{sor-R1}{\rm
See Section \ref{S:further} for a substitute for Assumption \ref{A2}
when transience is no longer assumed.
}
\end{remark}

\begin{remark}\label{sor-R3}{\rm
We will see in the next section that Assumption \ref{A2} implies
$E(x,r)$ and $E(y,r)$ are comparable if $d(x,y)\approx r$, but gives
no useful bounds when $d(x,y)\gg r$.
}
\end{remark}

Given $f\in \sF_\sG$ and $A\subset \sV$, define
\bee\label{sor-EA}
\sE_{\sG,A}=\tfrac12 \sum_{x,y\in A}[f(y)-f(x)]^2 C_{xy},
\eee
the Dirichlet form restricted to $A$.
Set $$f_A=\frac{1}{\mu(A)} \sum_{x\in A} f(x)\mu_x.$$

We say the adjusted Poincar\'e inequality (API) holds for $\sG$ if
there exists $\kappa_1>0$ and $\kappa_2>1$ such that
\bee\label{sor-API}
\sum_{y\in B(x,r)} [f(y)-f_{B(x,r)}]^2 \mu_y
\le \kappa_1 E(x,r)\sE_{\sG,B(x,\kappa_2r)}(f,f)
\eee
whenever $f\in L^2(\sV,\mu)$, $x\in \sV$, and $r\ge 1$.

\begin{remark}\label{sor-R4}{\rm
When $\sV=\Z^d$ with $\mu$ being counting measure
and $d\ge 3$, $V(x,r)\approx r^d$,
$C(x,r)\approx r^{d-2}$, and $E(x,r)\approx r^2$,
and we get the usual Poincar\'e inequality.
For a large class of nested fractals, $V(x,r)\approx r^{d_f}$,
$C(x,r)\approx r^{d_f-d_w}$, and $E(x,r)\approx r^{d_w}$,
where $d_f$ and $d_w$ are the fractal and walk dimensions, resp.
}
\end{remark}

We say the cut-off inequality (COI) holds for $\sG$ if there exist 
$\kappa_3, \kappa_4$, and $\th$ such that for each $x_0\in \sV$ and $R\ge 1$
there exists a function $\vp=\vp_{x_0,R}$ with the following properties.

(1) $\vp(x)\ge 1$ for $x\in B(x_0,R/2)$ and $\vp(x)=0$ for $x\notin B(x_0,R)$.

(2) For each $x,y\in \sV$,
$$|\vp(x)-\vp(y)|\leq \kappa_3\Big(\frac{d(x,y)}{R}\Big)^\th.$$

(3) If $1\le s\le R$ and $z\in \sV$, then
\begin{align}
\sum_{x\in B(z,s)} f(x)^2& \sum_y |\vp(y)-\vp(x)|^2 C_{xy}\label{sor-COI}\\
&\le \kappa_4\Big(\frac{s}{R}\Big)^{2\th}
\Big(\sE_{\sG, B(z,2s)}(f,f)+
E(z,s)^{-1} \sum_{x\in B(z,2s)} f(x)^2\mu_x\Big).\nn
\end{align}

\begin{remark}\label{sor-R5}{\rm
The COI is very similar to the CS inequality of \cite{BB}, where
an extensive  discussion can be found.
}
\end{remark}

We say a function $h$ on a subset $D$ of $\sV$ is harmonic if
$$\sL h(x)=0, \qq x\in D.$$
This is equivalent to
$$h(x)=\sum_y h(y) C_{xy}, \qq x\in D.$$
The elliptic Harnack inequality (EHI) holds for the weighted graph $\sG$ with
conductances $\{C_{xy}\}$ if there exists $c_1$
such that whenever $x_0\in \sV$, $r\ge 1$, and $h$ is non-negative and
harmonic in $B(x_0,2r)$, then
\bee\label{def-harm}
h(x)\le c_1h(y), \qq x,y\in B(x_0,r).
\eee

Our first main theorem is the following.

\bet\label{T1}
Suppose \eqref{p0} and Assumptions \ref{A1},  \ref{A2}, and \ref{A3} hold. 

(a) If
the EHI holds for $\sG$, then  both the API and COI hold for $\sG$.

(b)  If the API and COI hold
for $\sG$, then the EHI holds for $\sG$.
\eet

Suppose we have another set of conductances $\{C'_{xy}\}$ on the graph $\sG$.
We say $(\sG,C_{xy})$ and $(\sG, C'_{xy})$ are equivalent weighted graphs
if there exists $c_1<1$ such that
$$c_1 C_{xy}\le C'_{xy}\le c_1 C_{xy}, \qq x,y\in \sV.$$

Our second main theorem is the stability theorem.

\bet\label{T2}
Suppose $(\sG,C_{xy})$ and $(\sG, C'_{xy})$ are equivalent weighted
graphs.  
Suppose \eqref{p0} and   Assumptions \ref{A1}, \ref{A2}, and \ref{A3}
hold for $(\sG, C_{xy})$ and for $(\sG, C'_{xy})$.
If the EHI holds for $(\sG, C_{xy})$, 
then the EHI holds for 
 $(\sG,C'_{xy})$.
\eet

See Section \ref{S:further}
for a statement of these theorems in the context
of metric measure spaces or when Assumption \ref{A1} does not hold.

\section{Preliminaries}\label{S:prelim}

We introduce the cable process. Let $\sC$ consist of $\sV$ together
with copies of $(0,1)$, one for each edge in $\sG$. If $x\sim y$, we write
$(x,y)$ for the corresponding copy, and we call $(x,y)$ the cable
connecting $x$ and $y$. We identify $x$ with $0$ and $y$ with 1
on the cable connecting $x$ and $y$. We define $\mu(dz)$ be setting it equal to
$C_{xy}\, dz$ on the cable connecting $x$ and $y$, where $dz$ is linear
Lebesgue measure. If $x$ and $y$ are two points on the same cable or
one lies on a cable and the other is an endpoint of that cable, then we define
the distance between $x$ and $y$ by $|x-y|$. In $x$ and $y$ are on different
cables, we use $\min\{|x-z_x|+d(z_x,z_y)+|z_y-y|\}$ for the distance, where the minimum 
is taken over all vertices $z_x, z_y\in\sV$ such that $x$ in on a cable with
one end at $z_x$ and $y$ is on a cable with one end at $z_y$.
We continue to use the notation $d(x,y)$ for the distance and set
$$B'(x,r)=\{y\in \sC: d(x,y)<r\}, \qq V'(x,r)=\mu(B'(x,r)).$$

The cable process is the process that behaves like one-dimensional
Brownian motion speeded up deterministically by the factor $C_{xy}$ on
$(x,y)$ and when at a vertex $x$, picks the cable along which the next excursion
takes place according to the probabilities $C_{xy}/\mu_x$. More precisely,
if $x\in \sC-\sV$ and $x$ lies on  the cable $(y_0,y_1)$, let
$$\grad f(x)=\lim_{z\to x} \frac{f(z)-f(x)}{d(y_0,z)-d(y_0,x)}.$$
If $x\in \sV$ and $x\sim y$, let
$$\grad_y f(x)=\lim_{z\to x, z\in (x,y)} \frac{f(z)-f(x)}{d(x,z)}.$$
Since we only work with $|\grad f|$ and $|\grad_y f|$, we do not
need to be concerned with whether we use $y_0$ or $y_1$ in the definition
of $\grad f(x)$. Let
$$\sE_\sC(f,f)=\tfrac12 \int_{\sC-\sV} |\grad f(z)|^2\, \mu(dz),$$
let $\sF^0_\sC$ be the collection of continuous functions with compact
support such that $\grad f(z)$ exists at every point of $\sC-\sV$, $\grad_y f(x)$
exists at every $x\in \sV$ for which  $y\sim x$, and $|\grad f|$ is bounded. For
the domain of $\sE_\sC$, we use $\sF_\sC$, which is the completion of $\sF^0_\sC$
with respect to the norm
$$\Big(\int_\sC |f(z)|^2\, \mu(dz)\Big)^{1/2}+\sE_\sC(f,f)^{1/2}.$$
The cable process is the symmetric continuous Markov process $(Y_t, \P^x)$
corresponding
to $(\sE_\sC, \sF_\sC)$. Typically when constructing a process via Dirichlet
forms, there is a null set involved, and one has to talk about properties
holding quasi-everywhere. However, in our case $\P^x(Y_t \mbox{ ever hits } y)>0$
for each $x$ and $y$, and no null set is necessary.

Let $\sL_\sC$ be the infinitesimal generator of $Y$. See \cite{BB} for
a detailed description of $\sL_\sC$ and its domain.

\bep\label{prelim-P1}
Suppose \eqref{p0} and Assumptions \ref{A1} and \ref{A2} hold. Then 
$(\sE_\sC, \sF_\sC)$ is transient. Let 
$$C'(x,r)=\inf\{\sE_\sC(f,f): f|_{B'(x,r)}=1, f\in \sF_\sC\}$$
 be the capacity of 
$B'(x,r)$ and let $$E'(x,r)=\frac{V'(x,r)}{C'(x,r)}.$$ Then
there exist $c_1>0$, $\rho\in (0,1)$ and a positive integer $M$
such that
\begin{align}
V'(x,2r)&\le c_1 V'(x,r), \label{CVD}\\
C'(x,r)&\le \rho C'(x,2r), \label{CCG}\\
E'(x,r)&\le \rho E'(x,2r) \label{CEG}
\end{align}
whenever $x\in\sC$ and $r>0$.
Moreover there exists $M$ not depending on $x$ or $r$ such that
the boundary of $B'(x,r)$ can be covered by at most $M$ balls of radius
$r/8$.
\eep

\proof This follows easily by using the techniques of \cite[Section 3]{BB}
and we leave the details to the reader.
\qed

Given $f\in \sF_\sC$ and $A\subset \sC$, define
\bee\label{prelim-EA}
\sE_{\sC,A}=\tfrac12 \int_{A-\sV} |\grad f(x)|^2 \, \mu(dx).
\eee
and set $$f_A=\frac{1}{\mu(A)} \int_A f(x)\, \mu(dx).$$

We say the adjusted Poincar\'e inequality (API) holds for $\sC$ if
there exists $\kappa_1>0, \kappa_2>1$ such that
\bee\label{prelim-API}
\int_{B'(x,r)} [f(y)-f_{B'(x,r)}]^2 \, \mu(dy)
\le \kappa_1 E'(x,r)\sE_{\sC,B'(x,\kappa_2r)}(f,f)
\eee
whenever $f\in \sF_\sC$, $x\in \sC$.

We say the cut-off inequality (COI) holds for $\sC$ if there exist 
$\kappa_3, \kappa_4$, and $\th$ such that for each $x_0\in \sC$ and $R> 0$
there exists a function $\vp=\vp_{x_0,R}$ with the following properties.

(1) $\vp(x)\ge 1$ for $x\in B'(x_0,R/2)$ and $\vp(x)=0$ for $x\notin B'(x_0,R)$.

(2) For each $x,y\in \sC$,
$$|\vp(x)-\vp(y)|\leq \kappa_3\Big(\frac{d(x,y)}{R}\Big)^\th.$$

(3) If $0\le s\le R$ and $z\in \sC$, then
\begin{align}
\int_{ B'(z,s)}& f(x)^2  |\grad \vp(x)|^2 \, \mu(dx)\label{prelim-COI}\\
&\le \kappa_4\Big(\frac{s}{R}\Big)^{2\th}
\Big(\sE_{\sC, B'(z,2s)}(f,f)+
E'(z,s)^{-1} \int_{ B'(z,2s)} f(x)^2\, \mu(dx)\Big).\nn
\end{align}

We say a function $h$ in the domain of $\sL_\sC$  
is harmonic on a subset $D$ of $\sC$ if
$$\sL_\sC h(x)=0, \qq x\in D.$$
The elliptic Harnack inequality (EHI) holds for $\sC$ 
if there exists $c_1$
such that whenever $x_0\in \sC$, $r>0$, and $h$ is non-negative and
harmonic in $B'(x_0,2r)$, then
$$h(x)\le c_1h(y), \qq x,y\in B'(x_0,r).$$

\bep\label{prelim-P2}
(a) The COI holds for $\sC$ if and only the COI holds for $\sG$.

(b) The API holds for $\sC$ if and only the API holds for $\sG$.

(c) The EHI holds for $\sC$ if and only the EHI holds for $\sG$.
\eep

\proof
The proof of (a) is almost identical to that of Propositions 3.3 and 3.4 
of \cite{BB}. The same techniques
can be used to prove (b). (c) is \cite[Cor.~2.5]{BB}.
\qed

The main work in this paper is to prove the following.

\bet\label{T3}
Suppose \eqref{p0} and Assumptions \ref{A1}, \ref{A2}, and \ref{A3} hold. 

(a) If
the EHI holds for $\sC$, then  both the API and COI hold for $\sC$.

(b)  If both the API and COI
hold for $\sC$, then the EHI holds for $\sC$.
\eet

It will be clear from the context whether we are working with $\sC$ or $\sG$,
so henceforth we will drop the primes and write
$B(x,r), V(x,r), C(x,r)$, and $E(x,r)$ in place
of
$B'(x,r), V'(x,r), C'(x,r)$, and $E'(x,r)$, resp.
We write $\del B(x,r)$ for the boundary of $B(x,r)$.

\bel\label{prelim-L0}
There exists $c_1>0$ and $\rho'\in (0,1)$
such that
volume growth holds:
\bee\label{CVG}
V(x,r)\le \rho' V(x,2r), \qq x\in \sC, r>0;
\eee
capacity doubling holds:
\bee\label{CCD}
C(x,2r)\le c_1 C(x,r), \qq x\in \sC, r>0;
\eee
and expected occupation time doubling holds:
\bee\label{CED}
E(x,2r)\le c_1 E(x,r), \qq x\in \sC, r>0,
\eee
\eel

\proof Multiplying \eqref{CCG} and \eqref{CEG} together  
 gives volume growth.
Expected occupation time growth implies
$$C(x,2r)\leq \rho C(x,r)\frac{V(x,2r)}{V(x,r)},$$
and an application of volume doubling implies capacity
doubling. Finally, since $C(x,r)\le C(x,2r)$, volume doubling implies
$$\frac{V(x,2r)}{C(x,2r)}\leq c_2 \frac{V(x,r)}{C(x,r)},$$
which is expected occupation time 
doubling.
\qed

\bel\label{prelim-L1}
Let $a>0$. There exists $c_1$ depending on $a$ but not on $r$, $x$, or $y$
such that if $d(x,y)<ar$, then 
$$V(x,r)\le c_1 V(y,r), \quad C(x,r)\le c_1 C(y,r), \quad
E(x,r)\le c_1 E(y,r).$$
\eel

\proof
Since $B(x,r)\subset B(y, (1+a)r)$, volume doubling tells
us
$$V(x,r)\le V(y, (1+a)r)\le c_2 V(y,r),$$
and similarly for $V$ replaced by $C$. By symmetry, $C(y,r)\le c_2 C(x,r)$,
so taking the ratio, $E(x,r)\le c_2^2 E(y,r).$
\qed

In Proposition \ref{prelim-P1}
 we may without loss of generality assume that the center of 
each of the $M$ balls is within $r/8$ of $\del B(x,r)$. If we let $B_1, \ldots,
B_M$ be balls with the same centers but radii equal to $r/4$, then for
each $j\ge 2$, there exists $i<j$ and a point $y_j$ such that $y_j\in B_i\cap B_j$. 
If $h$ is non-negative and harmonic in $B(x,2r)-B(x, r/2)$ and the EHI
holds, then 
for $w\in B_i$ and $z\in B_j$,
$$h(w)\leq c_1 h(y_j)\le c_1^2 h(z).$$
Using this inequality at most $M$ times, there is thus a constant
$c_2$ such that if 
$y,z\in \del B(x,r)$, then
\bee\label{EHC}
h(y)\le c_2 h(z).
\eee

Let $G(x,y)$ be the Green function for the process $Y_t$. The
existence of $G$ is an easy consequence of Assumption \ref{A1}
and the structure of $\sC$. For $x$ fixed,
$h(z)=G(x,z)$ is a non-negative function that is harmonic in
$B(x,2r)-B(x,r/2)$ and so we may apply \eqref{EHC} to $G(x, \cdot)$
and obtain
\bee\label{EHC2}
G(x,y)\leq c_2G(x,z), \qq y,z\in \del B(x_0,r).
\eee

When the EHI holds, harmonic functions are H\"older continuous (see \cite{Moser1}),
and so there exist $c_3$ and $\beta$ such that if $h$ is harmonic in
$B(x_0,2r)$, then
\bee\label{HC}
|h(x)-h(y)|\le c_3 \Big(\frac{d(x,y)}{r}\Big)^\beta \Big(\sup_{B(x_0,2r)}
|h|\Big), \qq x,y\in B(x_0,r).
\eee

\section{Some consequences of the EHI}\label{S:conseq}

In this section we assume the EHI holds for a process $Y$ associated with
a Dirichlet form $(\sE, \sF)$.

The first estimate is standard.
Let $G(x,y)$ be the Green function for $Y$.

\bep\label{PE-P1}
There exists constants $c_1$ and $c_2$ such that if $r=d(x,y)$, then 
$$\frac{c_1}{C(x,r)}\leq G(x,y)\leq \frac{c_2}{C(x,r)}, \qq x\in \sC, r>0.$$
\eep

\proof Let $x$ and $y$ be fixed and let $r=d(x,y)$. Let $\nu$ be the
capacitary measure for $B(x,r)$. Then we know $\nu$ is supported on
$\del B(x,r)$, its total mass is $C(x,r)$, and $G\nu$ equals 1 on $B(x,r)$. 
(See \cite[Section II.5]{PTA}, for example. The proofs there are for Brownian motion but
are valid for any symmetric continuous strong Markov process.)
Using \eqref{EHC2},   we may write
\begin{align*}
1&=G\nu(x)=\int_{\del B(x,r)} G(x,z)\, \nu(dz) \ge
c_3G(x,y) \int_{\del B(x,r)} \, \nu(dz)\\
&=c_3 G(x,y) C(x,r).
\end{align*}
Rearranging gives the right hand inequality.
The left hand inequality is proved in the same way, replacing ``$\ge$''
by ``$\le$.''
\qed

Next we obtain an estimate on the time spent in $B(x,r)$. 

\bep\label{PE-P3}
There exist constants $c_1$ and $c_2$ such that
$$c_1 E(x,r)\le \int_{B(x,r)} G(x,z)\, \mu(dz)\leq c_2 E(x,r).$$
\eep

\proof 
Let $\rho'$ be the constant in Lemma \ref{prelim-L0}.
Applying \eqref{EHC2}, Proposition \ref{PE-P1}, and \eqref{CVG},
\begin{align*}
\int_{B(x,r)} G(x,z)\, \mu(dz)&\ge \int_{B(x,r)-B(x,r/2)} G(x,z)\, \mu(dz)\\
&\ge \frac{c_3}{C(x,r)} ( V(x,r)-V(x,r/2))\\
&\ge \frac{c_3(1-\rho')}{C(x,r)} V(x,r)\\
&=c_4 E(x,r).
\end{align*}
This gives the left hand inequality.

Similarly, we have
$$\int_{B(x,r)-B(x,r/2)} G(x,z)\, \mu(dz)\leq c_5 \frac{V(x,r)-V(x,r/2)}{C(x,r)}\le c_5 E(x,r)$$
for each $r>0$. We apply this with $r$ replaced by $2^{-k}r$ for $k=0,1, 
\ldots$, and sum. Using the fact that $Y$ spends 0 time at $x$ (locally $Y$ behaves
like a deterministic time change of Brownian motion), we obtain
\bee\label{G0}
\int_{B(x,r)} G(x,z)\, \mu(dz)
\le c_5 \sum_{k=0}^\infty E(x,2^{-k}r).
\eee
Using \eqref{CEG} repeatedly, we have $E(x, 2^{-k}r)\leq \rho^k E(x,r)$, so
$$\int_{B(x,r)} G(x,z)\, \mu(dz)\le c_5 E(x,r)\sum_{k=0}^\infty \rho^k,$$
which implies the right hand inequality.
\qed

\section{The adjusted Poincar\'e inequality}\label{S:api}

Let $G_D$ denote the Green function for $Y$ killed on exiting a domain
$D$. 

\bep\label{EiA-P1}
Suppose \eqref{p0}, Assumptions \ref{A1}, \ref{A2}, and \ref{A3}, and
the EHI hold.
There exists $k_0\ge 2$ and $c_1$ not depending on $x_0$ or $r$ such that if $r>0$
and $x,y\in B(x_0,r)$, then
$$G_{B(x_0, 2^{k_0}r)}(x,y)\geq \frac{c_1}{C(x_0,r)}.$$
\eep

\proof Let $s=d(x,y)$ and note $B(x,s)\subset B(x_0,4r)$.
By Proposition \ref{PE-P1} and \eqref{CCD}, there exists a constant $c_2$ such
that 
\bee\label{EiA-P1E1}
G(x,y)\ge \frac{c_2}{C(x,s)}\ge \frac{c_2}{C(x_0,4r)}
\ge \frac{c_3}{C(x_0,r)}.
\eee
By the strong Markov property, 
\bee\label{EiA-P1E2}
G_D(x,y)=G(x,y)-\E^x G(Y_{\tau_D},y),
\eee
where $\tau_D$ is the first time that $Y$ exits $D$. By \eqref{EHC2},
if $D=B(x_0,2^kr)$ for some $k\ge 1$ and 
$w\in \del D$, then
\bee\label{EiA-P1E3}
G(w,y)\le c_4G(w,x_0)\leq \frac{c_5 }{C(x_0,2^kr)}\leq
\frac{c_5\rho^k}{C(x_0,r)},
\eee
where $\rho$ is the constant in Proposition \ref{prelim-P1}. 
If we choose $k_0\ge 2$ large enough so that $c_5\rho^{k_0}\le c_3/2$
and combine \eqref{EiA-P1E1}, \eqref{EiA-P1E2}, and
\eqref{EiA-P1E3},
we then have
our proposition with $c_1=c_3/2$.
\qed

We write $(G_D)^2f$ for $G_D(G_Df)$.

\bep\label{EiA-P2} 
Suppose \eqref{p0}, Assumptions \ref{A1}, \ref{A2}, and \ref{A3}, and
the EHI hold.
Let $k_0$ be defined as in Proposition \ref{EiA-P1}
and let $D\allowbreak =B(x_0,2^{k_0}r)$. There exists $c_1$ not depending on 
$x_0$ or $r$ such that
$$(G_D)^2(x,y)\leq c_1E(x_0,r) G_D(x,y)$$
for all $x,y\in B(x_0,r)$.
\eep

\proof Write
$$(G_D)^2(x,y)=\int G_D(x,z) G_D(z,y)\, \mu(dz).$$
We let $s=d(x,y)$ (so that $s<2r$) and break the integral on the right
into integrals over $B(x,s/2)$ and over $B(x,s/2)^c$. 

For $z\in B(x,s/2)$, we have $d(z,y)\ge s/2$, and by \eqref{EHC2}
$$G_D(z,y)\leq c_2 G_D(x,y).$$
Since  $D\subset B(x, 2^{k_0+1}r)$, using 
Proposition \ref{PE-P3}, \eqref{CED}, and Lemma \ref{prelim-L1}
yields
\begin{align*}
\int_{B(x,s/2)} G_D(x,z)G_D(z,y)\, \mu(dz)
&\le c_2 G_D(x,y)\int_D G_D(x,z)\, \mu(dz)\\
&\le c_2 G_D(x,y) \int_D G(x,z)\,  \mu(dz)\\
&\le c_3 G_D(x,y)\int_{B(x, 2^{k_0+1}r)} G(x,z)\, \mu(dz)\\
&\le c_4 G_D(x,y) E(x, 2^{k_0+1}r)\\
&\le c_5 G_D(x,y) E(x,r)\\
&\le c_6 G_D(x,y) E(x_0,r).
\end{align*}

For $z\in B(x, s/2)^c$, we have $d(z,x)\ge s/2$, and by \eqref{EHC2}
$$G_D(x,z)\le c_2 G_D(x,y).$$
As above, using that $G_D$ is zero on $D^c$,
\begin{align*}
\int_{B(x,s/2)^c} G_D(x,z) G_D(z,y)\, \mu(dz)
&\le c_2 G_D(x,y) \int_D G_D(y,z)\, \mu(dz)\\
&\le  c_2 G_D(x,y)\int_D G(y,z)\, \mu(dz)\\
&\le c_2 G_D(x,y)\int_{B(y, 2^{k_0+1}r)} G(y,z)\, \mu(dz)\\
&\le c_7 G_D(x,y) E(y, 2^{k_0+1}r)\\
&\le c_8 G_D(x,y) E(y, r)\\
&\le c_9 G_D(x,y) E(x_0,r).
\end{align*}
In the third inequality we used the fact that $D\subset B(y, 2^{k_0+1}r)$,
and we used Lemma \ref{prelim-L1} for the last inequality.
Adding the integrals over $B(x,s/2)$ and $B(x,s/2)^c$ yields our result.
\qed

Let $G^\al$ be the $\al$-resolvent for $Y$  and $G^\al_D$ the $\al$-resolvent
for the process killed on exiting $D$.

\bep\label{EiA-P3} 
Suppose \eqref{p0}, Assumptions \ref{A1}, \ref{A2}, and \ref{A3}, and
the EHI hold.
Let $D$ be as in Proposition \ref{EiA-P2}. There exist  $c_1, c_2$ not depending on $x_0$ or $r$
such that if $\al=c_1/ E(x_0,r)$ and $x,y\in B(x_0,r)$, then
$$G_D^\al(x,y)\ge c_2/C(x_0,r).$$
\eep

\proof By the resolvent equation, $G^\al_D=G_D-\al G_DG_D^\al$, and so
$$G^\al_D(x,y)=G_D(x,y)-\al G_DG_D^\al(x,y)\ge
G_D(x,y)-\al (G_D)^2(x,y).$$
From Proposition \ref{EiA-P2} we know
$$(G_D)^2(x,y)\le c_3E(x_0,r) G_D(x,y)$$
for $x,y\in B(x_0,r)$.
By Proposition \ref{EiA-P1} we also know $G_D(x,y)\ge c_4/C(x_0,r)$.  
Then
\begin{align*}
 G^\al_D(x,y)&\ge G_D(x,y)(1-\al c_3 E(x_0,r))\\
&\ge \frac{c_4}{C(x_0,r)}(1-\al c_3 E(x_0,r)).
\end{align*}
If we take $c_1=(2c_3)^{-1}$, then
since $\al=c_1/E(x_0,r)$, we have
$1-\al c_3 E(x_0,r)\ge \frac12$, and our
result follows.
\qed

Given a ball $D$, we let $Y^r$ be the process $Y$ reflected on
the boundary of $D$. Since $Y$ behaves locally like a Brownian 
motion, it is clear how $Y^r$ can be described probabilistically.
Using a more analytic approach, $Y^r$ is the continuous symmetric
strong Markov process corresponding to $\sE_D$ with domain
$\{f\in \sF: \int_D (|f|^2 +|\grad f|^2)<\infty\}$.

\bet\label{EiA-T1}
Suppose \eqref{p0} and Assumptions \ref{A1}, \ref{A2}, and \ref{A3}  hold. If
the EHI holds for $\sC$, then the API holds for $\sC$.
\eet

\proof Fix $x_0$ and $r>0$. Let $B=B(x_0,r)$ and $D=B(x_0,2^{k_0}r)$,
where $k_0$ is as in Proposition \ref{EiA-P1}. 
Let $\al$ be as in Proposition \ref{EiA-P3}. Let $Y^r$ be the process $Y$ reflected on
 the boundary of $D$,  and let $G^\al_r$ be the $\al$-resolvent
for $Y^r$. Fix $f\in L^2(D)\cap \sF$. 
Take $x\in B$.
Then
\bee\label{EiA-E1}
\int_B (f(y)-f_B)^2\, \mu(dy)\leq \int_B (f(y)-\al G^\al_r f(x))^2\, \mu(dy).
\eee 
We have  for $x,y\in B$,
\begin{align*}
\al G^\al_r (x,y)&\ge \al G^\al_D(x,y)\geq \frac{c_1}{E(x_0, r)C(x_0, r)}\\
&\ge \frac{c_1} {V(x_0, r)}.
\end{align*}
For any function $h$,
$$\al G^\al_rh(x)=\int_D h(y) \al G^\al_r(x,y)\, \mu(dy)\ge \frac{c_1}{V(x_0,r)}
\int_B h(y)\, \mu(dy),$$
and letting $h(y)=(f(y)-f_B)^2$, we obtain
\begin{align}
\int_B (f(y)-f_B)^2\, \mu(dy)&
\le c_2V(x_0,r) [\al G^\al_r((f(\cdot)-\al G^\al_r f(x))^2)(x)]\label{EiA-E1.5}\\
&=c_2 V(x_0,r) [\al G^\al_r (f^2)(x)-(\al G^\al_r f(x))^2].\nn
\end{align} 
The right hand side is non-negative. 
Integrating both sides over the set $D$ with respect to
the measure $\mu(dx)$,  multiplying by $\mu(B)^{-1}$, and
using volume doubling gives
\begin{align}
\int_B&(f(y)-f_B)^2\, \mu(dy)
\label{EiA-E2}\\
&\le 
c_2\Big[\int_D \al G^\al_r(f)^2(x)\, \mu(dx)-\int_D (\al G^\al_r f(x))^2\, \mu(dx)
\Big].\nn
\end{align}
If $\angel{\cdot,\cdot}$ is the inner product with respect to $L^2(D)$, then
using the symmetry of the resolvent, the first integral inside the brackets on the last line is
$$\angel{\al G^\al_r(f^2),1}=\angel{f^2, \al G^\al_r 1}=\angel{f^2,1}=
\norm{f}^2_{2},$$
where we write $\norm{\cdot}_2$ for the $L^2$ norm on $D$.  
The second integral on the last line of \eqref{EiA-E2}  is $\norm{\al G^\al_r f}^2_{2}$,
and  we thus have
\bee\label{EiA-E3}
\int_B(f(y)-f_B)^2\, \mu(dy)\leq c_2[\, \norm{f}^2_2-\norm{\al G^\al_r f}^2_2].
\eee

We now use the spectral theorem for $L^2(D)$. Let $\{E_\lam\}$ be
the spectral resolution of the operator $\sL^r$, the infinitesimal
generator of $Y^r$. Each $E_\lam$ is a projection, and we can write
$$f=\int_0^\infty \, dE_\lam f, \qq \norm{f}_2^2=\int_0^\infty
\, d\angel{E_\lam f, E_\lam f}.$$
For $f\in \sF$, we have
$$\sE_D(f,f)=\int_0^\infty \lam \, d\angel{E_\lam f, E_\lam f}.$$
We also  have
$$\al G^\al_r f=\int_0^\infty \frac{\al}{\al+\lam}\, dE_\lam f, 
\qq \norm{\al G^\al_r f}_2^2=\int_0^\infty \Big(\frac{\al}{\al+\lam}\Big)^2
\, d\angel{E_\lam f, E_\lam f}.$$
Since
$$1-\Big(\frac{\al}{\al+\lam}\Big)^2=\frac{2\lam(\al+\lam/2)}{(\al+\lam)^2}
\leq \frac{2\lam}{\al},$$
then 
\begin{align}
\norm{f}_2^2-\norm{\al G^\al_rf}_2^2&=\int_0^\infty \Big(1-\Big(\frac{\al}{\al+\lam}\Big)^2\Big)
\, d\angel{E_\lam f, E_\lam f}\label{EiA-56}\\
&\le c_2\frac{2}{\al} 
\int_0^\infty \lam \, d\angel{E_\lam f, E_\lam f}\nn\\
&= c_3 E(x_0,r) \sE_D(f,f).\nn
\end{align}
Combining \eqref{EiA-E3} and \eqref{EiA-56} proves the API.
\qed

\section{Proofs of main theorems}\label{S:proofs}

Throughout we assume \eqref{p0} and Assumptions \ref{A1}, \ref{A2},
and \ref{A3}.
We continue the cable system context unless stated otherwise.

We need two propositions which will be used to show that the COI and API
imply the EHI.

Fix $x_0\in \sC$, let $R\ge 1$, and let $\vp$ be the cut-off function given
by the COI. Let
$$\gamma=1+ E(x_0,R) |\grad \vp|^2.$$

\bep\label{HI-P1} Suppose the API holds for $\sC$  with constants 
$\kappa_1$ and $\kappa_2$ and also the COI holds. 
Let $x\in B(x_0,R)$, let $I=B(x,s)$ with $s\le R$, and 
let $I^*=B(x,2s)$, $I^{**}=B(x,2\kappa_2s)$. Suppose $f$ and its gradient are square integrable over $I^{**}$ and let
$f_A=\mu(A)^{-1}\int _A f\, d\mu$.
Then
\bee\label{HI-E1}
\int_I f^2\gamma\leq c_1 (s/R)^{2\th} E(x_0,R)
\Big(\int_{I^*} |\grad f|^2+ E(x,s) ^{-1}\int_{I^*} f^2\Big)
\eee
and
\bee\label{HI-E2}
\int_I (f-f_{I^*})^2\gamma\leq c_2(s/R)^{2\th} E(x_0,R) \int_{I^{**}} |\grad f|^2.
\eee
If $J\subset I$, then 
$$\int_J f^2\gamma \leq c_3\Big(E(x_0,R)(s/R)^{2\th}\Big)
\int_{I^{**}} |\grad f|^2+\mu(J)^{-1}\Big(\int_J |f|\gamma\Big)^2.$$
Finally,
$$\int_{B(x_0,R)}\gamma \leq c_4 V(x_0,R).$$
\eep

\proof The condition \eqref{CEG}  implies that $E(x,R)/E(x,s)\geq c_5 (R/s)^\beta$ for some $\beta>0$ and $c_5>0$ not depending on $x,R$, or $s$.
Without loss of generality we may assume $2\th< \beta$.
Then $$(s/R)^{2\th} E(x,R) E(x,s)^{-1}\ge c_6$$
since $s\le R$. Using Lemma \ref{prelim-L1}, $E(x_0,R)\ge c_7 E(x,R)$ and hence
\begin{align*}
\int_I f^2 \gamma&=\int_I f^2 + E(x_0,R)\int _I f^2|\grad \vp|^2\\
&\leq \int_I f^2 +c_8(s/R)^{2\th} E(x_0,R)\int_{I^*} |\grad f|^2
+c_8 (s/R)^{2\th} \frac{E(x_0,R)}{E(x,s)}\int_{I^*} f^2\\
&\le c_9(s/R)^{2\th} E(x_0,R)\int_{I^*} |\grad f|^2
+c_9 (s/R)^{2\th} \frac{E(x_0,R)}{E(x,s)}\int_{I^*} f^2.
\end{align*} 
Applying this to $f-f_{I^*}$, we have
$$\int_I (f-f_{I^*})^2 \gamma
\leq c_{10} (s/R)^{2\th} E(x_0,R)
\Big(\int_{I^*} |\grad f|^2 + E(x,s)^{-1}\int_{I^*}(f-f_{I^*})^2\Big).$$
Applying the API to $B(x,2s)$,
$$E(x,s)^{-1}\int_{I^*}(f-f_{I^*})^2\le c_{11} \int_{I^{**}} |\grad f|^2.$$
Combining gives \eqref{HI-E2}.

The remainder of the proof is exactly as in \cite[Prop.~5.2]{BB}.
\qed
 
Here is a substitute for \cite[Prop.~5.7]{BB}.

\bep\label{HI-P2} Suppose the API holds for $\sC$  with constants $\kappa_1$
and $\kappa_2$ and also the COI holds. 
Let $S>0$ and let $u$ be positive and harmonic in $B(x_0,2\kappa_2S)$ and let $w=\log u$.
Then 
$$\int_{B(x_0,2S)} |\grad w|^2\, d\mu\le c_1 C(x_0,S).$$
\eep

\proof Let $\vp_1$ be the cut-off function for $B(x_0,2\kappa_2S)$ given by the COI.
Exactly as in the proof of \cite[Prop.~5.7]{BB} we have
$$\int_{B(x_0,2S)} |\grad w|^2\, d\mu\le 
\int \vp_1^2 |\grad w|^2 \, d\mu \le c_2\int |\grad \vp_1|^2\, d\mu.$$
Applying the COI in $B(x_0,2\kappa_2S)$ with $f=1$ and $s=2\kappa_2S$ yields
$$\int |\grad \vp_1|^2\le c_3 E(x_0,s)^{-1}\int_{B(x_0,2s)} \, d\mu
=c_3 V(x_0,4\kappa_2S)/E(x_0,2\kappa_2S).$$
Using \eqref{CVD} and \eqref{CCD}  yields our result.
\qed

Combining with \eqref{HI-E2} tells us that
\bee\label{HI-E4}
\int_{B(x_0,R)} |w-w_{B(x_0,R)}|^2 \gamma \le c_4 E(x_0,R) C(x_0,R)
=c_4 V(x_0,R).
\eee

\longproof{of Theorem \ref{T3}}
We proved that the EHI for $\sC$ implies the API for $\sC$ in Theorem \ref{EiA-T1}.

That the EHI for $\sC$  implies the COI  for $\sC$ is proved in almost the identical way that
it is done in \cite[Sect.~4]{BB}. We replace the use of $\psi(r)$
there by $E(x_0,r)$ and also replace appearances of $r^\beta$ by $E(x_0,r)$
The analogue of Lemma 4.7(a) of \cite{BB} follows from
Proposition \ref{PE-P1}. To prove the analogue of \cite[Lemma 4.7(b)]{BB},
we use Proposition \ref{EiA-P1} and then follow the proof given in \cite{BB}.

Away from the Green function is H\"older continuous in each variable 
by \eqref{HC}.  The FVG condition
of \cite{BB} is implied by our current volume growth condition.

With Propositions \ref{HI-P1} and \ref{HI-P2}  in place of 
Propositions 5.2 and 5.7  of 
\cite{BB}, 
we can follow the argument of \cite[Section 5]{BB} to show
that the API and COI together imply the EHI.
\qed

\longproof{of Theorem \ref{T1}}
If \eqref{p0}, Assumptions \ref{A1}, \ref{A2}, and 
\ref{A3}, and the EHI hold for $(\sG,C_{xy})$,
Propositions \ref{prelim-P1} and \ref{prelim-P2}
tell us that the corresponding facts hold for the cable system $\sC$.
By Theorem \ref{T3}, the API and COI hold for $\sC$, and by
Proposition \ref{prelim-P2} again, the API and COI hold for the
weighted graph. This proves (a).
The proof of (b) is similar.
\qed

\longproof{of Theorem \ref{T2}}
Suppose \eqref{p0} and  Assumptions \ref{A1}, \ref{A2}, and \ref{A3}
 hold for
$(\sG, C_{xy})$ and for $(\sG, C'_{xy})$. Suppose the EHI
holds for $(\sG, C_{xy})$. Then by Theorem \ref{T1} the API and COI hold for
$(\sG, C_{xy})$. Since
 $(\sG, C_{xy})$ and $(\sG, C'_{xy})$ are
equivalent weighted graphs, then capacities of balls are comparable, and
hence expected occupation times are comparable. Therefore the API and COI hold for $(\sG,C'_{xy})$.
By Theorem \ref{T1}, the EHI holds for $(\sG, C'_{xy})$.
\qed

\section{Further results}\label{S:further}

\subsection{The general case}

We now consider the general case for infinite graphs. Theorem \ref{T5}
also can be used in the transient case.

For $x\in \sV$ and $r\ge 1$, let $\wt C(x,r)$ be the capacity of $B(x,r)$
with respect to the process killed on exiting $B(x,8r)$. Thus
$$\wt C(x,r)=\inf\{\sE_\sG(f,f): f|_{B(x,r)}=1, f|_{B(x,8r)^c}=0, f\in \sF\}.$$
Let $\wt E(x,r)=V(x,r)/\wt C(x,r)$. We assume \eqref{p0}, volume doubling,
expected occupation time growth
(for $\wt E$), and that the boundary of $B(x,r)$ can be covered
by at most $M$ balls of radius $r/8$.
$\wt C(x,r)$ is no longer necessarily monotone in $r$, and so we must
make an additional assumption, that of capacity comparability:
there exists $c_1$  not depending on $x,y,$ or $r$ such that if $d(x,y)<2r$,
then $$c_1\wt C(x,r) \le \wt C(y,2r)\le c_1^{-1} \wt C(x,r).$$
In particular, taking $x=y$ shows that $\wt C(x,r)$ and $\wt C(x,2r)$
are comparable.
This implies expected occupation time comparability:
there exists $c_2$ such that 
\bee\label{DEC}
c_2\wt E(x,r) \le \wt E(y,2r)\le c_2^{-1} \wt E(x,r).
\eee

Now define the API and COI in terms of $\wt E$ instead of $E$.

\bet\label{T5}
Suppose  \eqref{p0}, Assumption \ref{A3}, volume doubling, expected
occupation time growth, and capacity comparability hold for $\sG$. 

(a) If the EHI holds, then the API and COI hold.

(b) If  the API and COI hold, then the EHI holds.

(c) Let  $(\sG, C_{xy})$ and $(\sG, C'_{xy})$ be equivalent graphs. 
 Suppose  \eqref{p0}, Assumption \ref{A3}, volume doubling, expected
occupation time growth, and capacity comparability also hold for $(\sG, C'_{xy})$.
If the EHI holds for $(\sG, C_{xy})$, then it holds for $(\sG, C'_{xy})$.
\eet

\proof As in the proofs of Theorems \ref{T1} and \ref{T2}, we immediately
transfer to the cable system.
The proof of Proposition \ref{PE-P1} still applies and we have that 
$G_{B(x,r)}(x,y)$ is comparable to  $1/C(x,r)$.
The proof of Proposition \ref{PE-P3} shows that $\int_{B(x,r)} G_{B(x,8r)}(x,z)
\, \mu(dz)$ is comparable to $\wt E(x,r)$. 

For $x_0\in \sC$ and $r>0$, let $D=B(x_0,8r)$. Then if
$x,y\in B(x_0,r)$, we have
$$(G_D)^2(x,y)\leq c_1\wt E(x_0,r) G_D(x,y).$$
The proof of this is the same as the proof of Proposition 
\ref{EiA-P2}, but we use \eqref{DEC} to compare $\wt E(x,r)$
and $\wt E(y,r)$.
We then conclude 
$$G_D^\al(x,y)\geq c_2/\wt C(x_0,r),$$
just as in the proof of Proposition \ref{EiA-P3}.
We then argue that the EHI implies the API as in the proof of
Theorem \ref{EiA-T1}.
The remainder of the proof of Theorem \ref{T5}
is as in Section \ref{S:proofs}.
\qed

\subsection{Metric measure spaces}

There is no difficulty extending our theorems 
to more general continuous state spaces. 
See \cite{BBK} for the definitions of all terms introduced
in this subsection. 
Let $(X,d,\mu)$ be a
metric measure space such that the metric is geodesic and $X$
has infinite diameter. 
Examples of such spaces include Riemannian manifolds,
cable systems, Euclidean domains with smooth boundary, and fractals. 

Let $(\sE, \sF)$ be a local regular
Dirichlet form. 
Associated to  $f\in \sF\cap L^\infty$ is a measure
$\Gamma(f,f)(dx)$ characterized by
$$\int _X \wt g(x) \, \Gamma(f,f)(dx)=2\sE(f,fg)-\sE(f^2,g)$$
for all $g\in \sF\cap L^\infty$, where
$\wt g$ is the quasi-continuous modification of $g$.
Define
$$\sE_A(f,f)=\int_A \, \Gamma(f,f)(dx).$$

Let $B(x,r)$ be the ball of radius $r$, $V(x,r)=\mu(B(x,r))$.
Assume $(\sE, \sF)$ is transient, let
$$C(x,r)=\inf\{\sE(f,f): f|_{B(x,r)}=1, f\in \sF\},$$
and $E(x,r)=V(x,r)/C(x,r)$.
Assume that Assumption \ref{A2} holds; the statement in the present
context is the same as the one in Section \ref{S:prelim} provided we drop
the primes.  Again dropping the primes, define the API, COI, and 
EHI as in Section \ref{S:prelim}. Assume the analogue of Assumption \ref{A3}.  
We need one more regularity condition, namely, that the associated
continuous symmetric strong Markov process spends 0 time at any given 
point, or equivalently, for each $x$, 
\bee\label{F-A}
G1_{B(x,r)}(x)\to 0 \mbox{ as } r\to 0,
\eee
where here  $G$ is the Green potential operator.

We then have the analogues of Theorems \ref{T1} and \ref{T2}. We say
two Dirichlet forms $\sE$ and $\sE'$ are equivalent if they have
the same domain $\sF$ and there exists $c_1$ such that
$$c_1\sE(f,f)\le \sE'(f,f)\le c_1^{-1}\sE(f,f), \qq f\in \sF.$$

\bet\label{T4}
Assume that the analogues of Assumptions \ref{A1}, \ref{A2}, and \ref{A3},
and \eqref{F-A} hold for $(\sE, \sF)$. 

(a) If the EHI holds, then the API
and COI hold for $(\sE,\sF)$.

(b) If the API and  COI hold for $(\sE, \sF)$, then
the EHI holds for $(\sE, \sF)$.

(c) Let $\sE$ and $\sE'$ be equivalent. Assume that the analogues of Assumptions \ref{A1}, \ref{A2}, and \ref{A3} and
\eqref{F-A} hold for $(\sE', \sF)$. 
If the EHI holds for $\sE$, then it holds for $\sE'$.
\eet

\proof We modify the proof of Theorem \ref{T3} in a manner entirely
similar to the way \cite{BBK} extended the results of \cite{BB} to
metric measure spaces. \eqref{F-A} comes in when deriving \eqref{G0}.
The details are left to the interested reader.
\qed

\begin{remark}\label{proofs-R1}{\rm
We can similarly state and prove the analogue of Theorem \ref{T5}.
}
\end{remark}


\medskip

\ni {\bf Richard F. Bass}\\
Department of Mathematics\\
University of Connecticut \\
Storrs, CT 06269-3009, USA\\
{\tt bass@math.uconn.edu}
\ms

\end{document}